\documentclass[12pt]{article} 
\usepackage{amsmath}
\usepackage{amsfonts}
\usepackage{enumerate}
\usepackage{fancyhdr}
\usepackage{amsthm}
\usepackage{bm}
\numberwithin{equation}{section}
\usepackage{typearea}
\typearea{15}

\newtheorem{theo}{Theorem}[section]
\newtheorem{defi}{Definition}[section]
\newtheorem{lemm}{Lemma}[section]

\newtheorem{rem}{Remark}[section]
\newtheorem{con}{Question}[section]
\newtheorem*{ack}{Acknowledgments}
\newtheorem*{pf}{Proof}

\title{Real zeros of the Barnes double zeta function in the interval $(1, 2)$}
\author{Kazuma Sakurai}
\date{}

\begin{document}
\maketitle
\abstract
Let $a, w_1, w_2,\cdot\cdot\cdot, w_r >0$  and $s \in \mathbb{C}$. We put $\bm{w} = (w_1,\cdot\cdot\cdot,w_r)$. Then the Barnes $r$-ple zeta function is defined by $\zeta_r(s, \bm{w}, a) = \sum_{m_1=0}^{\infty} \cdot\cdot\cdot \sum_{m_r=0}^{\infty} 1/(a+m_1w_1+\cdot\cdot\cdot +m_rw_r)^s$ when $\sigma := \Re(s)>r$. In this paper, we show that the Barnes double zeta function $\zeta_2(\sigma, \bm{w}, a)$ has real zeros in the interval $(1,2)$ if and only if $0< a < (w_1+w_2)/2$ and the number of such zero is precisely one if $0< a< (w_1+w_2)/2$.

\section{Introduction}

~~~~~The Hurwitz zeta function was first introduced by Hurwitz as a generalization of the Riemann zeta function $\zeta(s) = \sum_{n=1}^{\infty} n^{-s}$.
\begin{defi} {\rm (}see {\rm \cite[Chapter 9]{aik})}
	Let $0 < a \le 1$, $s \in \mathbb{C}$. Then the Hurwitz zeta function $\zeta(s, a)$ is defined by 
	\begin{equation*}
		\zeta(s, a) := \sum_{n=0}^{\infty} \frac{1}{(n+a)^s},~~~~~s := \sigma + i t, \sigma>1, t \in \mathbb{R}.
	\end{equation*}
\end{defi}
As regards the real zeros of the Hurwitz zeta function, Nakamura \cite{naka1, naka2} investigated them in the interval $(0, 1)$ and $(-1, 0)$, and gave the following conditions.
\begin{theo}\label{n1}{\rm (\cite{naka1, naka2})} Let $b_2^{\pm}:=(3 \pm \sqrt{3})/6$. Then
\begin{enumerate}[{\rm (1)}]
	\item $\zeta(\sigma, a)$ has zeros in $(0, 1)$ if and only if $0<a<1/2$.
	\item $\zeta(\sigma, a)$ has zeros in $(-1, 0)$ if and only if 
	$0<a<b_2^{-}$ or $1/2 < a < b_2^{+}$.
\end{enumerate}
\end{theo}
Furthermore, Matsusaka \cite{matsusaka} extended this result to general intervals in the negative real numbers.
\begin{theo}\label{m1}{\rm (\cite{matsusaka})}
	Let $N \ge -1$ be an integer. Then $\zeta(\sigma, a)$ has zeros in $(-N-1, -N)$ if and only if $B_{N+1}(a)B_{N+2}(a) <0$, where $B_n(a)$ is the nth Bernoulli polynomial.
\end{theo}
We can easily see that Theorem \ref{n1} is a special case of Theorem \ref{m1}. Furthermore, we can state the latter theorem explicitly as follows.
\begin{theo}\label{m2}{\rm (\cite{matsusaka})}
	Let $N \ge 0$ be an integer. Then $\zeta(\sigma, a)$ has zeros in $(-N-1, -N)$ if and only if 
\[
	\begin{cases}
		0<a<b_{N+2}^{-}~or~1/2<a<b_{N+2}^{+}~~~~if~N~is~even, \\
		b_{N+1}^{-}<a<1/2~or~b_{N+1}^{+}<a<1~~~~if~N~is~odd
	\end{cases}
\]
where $b_n^{\pm}$ are the two roots of $B_n(x)$ in $(0, 1)$.
\end{theo}
For $N \ge 4$, we can show the uniqueness of the zero in each interval $(-N-1, -N)$. More recently, Endo and Suzuki \cite{es} proved the uniqueness of the zero in (0, 1) and found its asymptotic behavior with respect to a.

\begin{theo}\label{ensu}{\rm (\cite{es})} We have the following:
	\begin{enumerate}[{\rm (1)}]
		\item The Hurwitz zeta-function $\zeta(\sigma, a)$ has precisely one zero in the interval $(0, 1)$ when $0<a<1/2$, and these zeros are all simple.
		\item For $0 < a < 1/2$, let $\beta(a)$ denote the unique zero of $\zeta(\sigma, a)$ in $(0, 1)$. Then $\beta : (0, 1/2) \rightarrow (0, 1)$ is  a strictly decreasing $C^{\infty}$-diffeomorphism. Moreover, we have an asymptotic expansion
		\[
			\beta(a) = 1-a+a^2 \log a + O(a^2)~~as~~a \rightarrow 0+.
		\]
	\end{enumerate}
\end{theo}

\begin{rem}
\begin{rm}
	More generally,  in \cite{tok} Tokarev has obtained a generalization of Descartes' rule of signs which provides an upper bound for the number of positive zeros of a function. \end{rm}
\end{rem}
We are interested in extending these results to the Barnes $r$-ple zeta function. First, we define the Barnes $r$-ple zeta function.
\begin{defi}\label{bz}{\rm (}see {\rm \cite[Chapter 13]{aik})}
	Let all $a, w_1,w_2,\cdot\cdot\cdot.w_r$ be positive real numbers. We put $\bm{w}=(w_1,w_2,\cdot\cdot\cdot,w_r)$. Then the Barnes $r$-ple zeta function $\zeta_r(s,\bm{w}, a)$ is defined by
	\begin{equation*}
		\zeta_r(s,\bm{w}, a):=\sum_{m_1=0}^{\infty} \cdot\cdot\cdot \sum_{m_r=0}^{\infty} \frac{1}{(a+m_1w_1+\cdot\cdot\cdot +m_rw_r)^s} , ~~s := \sigma + it,~\sigma > r,~t \in \mathbb{R}.
	\end{equation*}
\end{defi}
This series converges absolutely in the half-plane $\sigma > r$ and uniformly in each compact subset of this half-plane. Moreover $\zeta_r(s, {\bm w}, a)$ can be meromorphically continued to the whole complex plane with simple poles at $s = 1, 2, \cdot \cdot \cdot , r$. We easily see that the Barnes $r$-ple zeta function has no real zeros in the half-plane $\sigma > r$. In the particular case $r=1$, this is essentially the Hurwitz zeta function. Namely, we have 
\begin{equation*}
	\zeta_1(s, w_1, a) = w_1^{-s} \zeta(s, \frac{a}{w_1}).
\end{equation*}
In this case, we see from Theorem \ref{m1} that $\zeta_1(\sigma, w_1, a)$ has zeros in $(-N-1, -N)$ for $N \ge -1$ if and only if $B_{N+1}(a/w_1)B_{N+2}(a/w_1) <0$. 

In this paper, we extend (1) of Theorem \ref{n1} and (1) of Thenorem \ref{ensu} to the Barnes double zeta function. 

\begin{theo}\label{main}We have the following: 
\begin{enumerate}[{\rm (1)}]
	\item $\zeta_2(\sigma, {\bm w}, a)$ has zeros in $(1,2)$ if and only if $0 < a < (w_1+w_2)/2$.
	\item When $0 < a < (w_1 + w_2)/2$, $\zeta_2(\sigma, {\bm w}, a)$ has precisely one simple zero in $(1, 2)$.
\end{enumerate}
\end{theo}	

\section{Proof of Theorem \ref{main}}
In this section, we first show (1) of Theorem \ref{main} by using Nakamura's method (see \cite[{\S} 2.2]{naka1}).
In order to prove Theorem \ref{main}, we define $G_2(a, \bm{w}, t)$ by
	\begin{eqnarray*}
		G_2(a,\bm{w},t) &:=& \frac{e^{(w_1+w_2-a)t}}{(e^{w_1t}-1)(e^{w_2t}-1)}-\frac{1}{w_1w_2t^2}. 
	\end{eqnarray*}

\begin{defi}\label{Barnbern}{\rm (see \cite[p.~377]{barnes})} Define the generalized r-ple Bernoulli polynomial $B_k^{(r)}(x,\bm{w})$ $(k \in \mathbb{N}_0)$ by
\begin{eqnarray*}
	\frac{t^re^{(w_1+\cdot\cdot\cdot + w_r-x)t}}{(e^{w_1t}-1)\cdot\cdot\cdot (e^{w_rt}-1)}&=& \sum_{k=0}^{\infty}B_k^{(r)}(x,\bm{w}) \frac{t^k}{k!} \\
&=& \frac{1}{w_1\cdot\cdot\cdot w_r}+ \frac{-2x + w_1+\cdot\cdot\cdot+ w_r}{2w_1\cdot\cdot\cdot w_r}t + \cdot\cdot\cdot,
\end{eqnarray*}
whose expression is valid for $|t|< \min \left\{2\pi/w_1, ... , 2\pi/w_r \right\}$ . In particular $B_k^{(1)}(x, 1) = B_k(1-x) = (-1)^kB_k(x)$.
\end{defi}

\begin{lemm}\label{inte1} For $1<\sigma<2$ we have the integral expression
	\begin{eqnarray}\label{inte}
		\Gamma(s)\zeta_2(s,\bm{w},a) &=& \int_{0}^{\infty}\left( \frac{e^{(w_1+w_2-a)t}}{(e^{w_1t}-1)(e^{w_2t}-1)}-\frac{1}{w_1w_2t^2}\right)t^{s-1}dt \nonumber \\
		&=& \int_{0}^{\infty} G_2(a, \bm{w}, t)t^{s-1} dt. 
	\end{eqnarray}
\end{lemm}

\begin{pf}
\begin{rm}
Let  $0< \lambda < \min \left\{2\pi/w_1,  2\pi/w_2 \right\}$. When $\sigma >2$, it is well-known that 
\begin{equation}\label{積分表示}
	\zeta_2(s, {\bm w}, a) =\frac{1}{\Gamma(s)} \int_{0}^{\infty}  \frac{e^{(w_1+w_2-a)t}}{(e^{w_1t}-1)(e^{w_2t}-1)}t^{s-1}dt
\end{equation}
{\rm (see \cite[p. 210]{aik})}. Hence we have
\begin{eqnarray}\label{a}
	\Gamma(s)\zeta_2(s,\bm{w}, a)&=&\int_{0}^{\infty}  \frac{e^{(w_1+w_2-a)t}}{(e^{w_1t}-1) (e^{w_2t}-1)}t^{s-1}dt  \nonumber \\
	&=& \int_{0}^{\lambda}  \frac{e^{(w_1+w_2-a)t}}{(e^{w_1t}-1) (e^{w_2t}-1)}t^{s-1}dt +\int_{\lambda}^{\infty}  \frac{e^{(w_1+w_2-a)t}}{(e^{w_1t}-1) (e^{w_2t}-1)}t^{s-1}dt \nonumber \\
	&=& \int_{0}^{\lambda}  \left(\frac{e^{(w_1+w_2-a)t}}{(e^{w_1t}-1) (e^{w_2t}-1)}-\frac{1}{w_1w_2t^2}\right)t^{s-1}dt  \nonumber \\
	&&+\int_{\lambda}^{\infty}  \frac{e^{(w_1+w_2-a)t}}{(e^{w_1t}-1) (e^{w_2t}-1)}t^{s-1}dt +\frac{1}{w_1w_2}\int_{0}^{\lambda} t^{s-3}dt\nonumber \\
	&=&  \int_{0}^{\lambda}  \left(\frac{e^{(w_1+w_2-a)t}}{(e^{w_1t}-1) (e^{w_2t}-1)}-\frac{1}{w_1w_2t^2}\right)t^{s-1}dt  \nonumber \\
	&&+\int_{\lambda}^{\infty}  \frac{e^{(w_1+w_2-a)t}}{(e^{w_1t}-1) (e^{w_2t}-1)}t^{s-1}dt +\frac{1}{w_1w_2}\frac{\lambda^{s-2}}{s-2}.
\end{eqnarray}
By Definition \ref{Barnbern}, for $\sigma > 1$, it holds that
\begin{eqnarray*}
	&& \int_{0}^{\lambda}  \left|\frac{e^{(w_1+w_2-a)t}}{(e^{w_1t}-1) (e^{w_2t}-1)}-\frac{1}{w_1w_2t^2}\right||t^{s-1}|dt \\
	 &\ll& \int_{0}^{\lambda} t^{-1}t^{\sigma-1} dt.
\end{eqnarray*}
On the other hand, we have 
\begin{equation*}
	\frac{1}{w_1w_2}\frac{\lambda^{s-2}}{s-2} = -\int_{\lambda}^{\infty} \frac{1}{w_1w_2} \frac{1}{t^2} t^{s-1}dt,  ~~~~~1 <\sigma < 2.
\end{equation*}
Therefore, the integral expression
\begin{eqnarray*}
	\Gamma(s)\zeta_2(s,\bm{w}, a) &=& \int_{0}^{\lambda}  \left(\frac{e^{(w_1+w_2-a)t}}{(e^{w_1t}-1) (e^{w_2t}-1)}-\frac{1}{w_1w_2t^2}\right)t^{s-1}dt -\int_{\lambda}^{\infty} \frac{1}{w_1w_2} \frac{1}{t^2} t^{s-1}dt \nonumber \\
	&&+\int_{\lambda}^{\infty}  \frac{e^{(w_1+w_2-a)t}}{(e^{w_1t}-1) (e^{w_2t}-1)}t^{s-1}dt
\end{eqnarray*}
gives an analytic continuation for $1< \sigma< 2$. Thus we obtain this lemma.
\qed \end{rm}
\end{pf}
\begin{lemm}\label{inte2} Let  $0< \lambda < \min \left\{2\pi/w_1,  2\pi/w_2 \right\}$. For $\sigma>0$ we have the integral expression
	\begin{eqnarray*}
	\begin{split}
		\Gamma(s)\zeta_2(s,\bm{w},a)&= \int_{\lambda}^{\infty} \frac{e^{(w_1+w_2-a)t}}{(e^{w_1t}-1)(e^{w_2t}-1)}t^{s-1}dt \\
		&\quad + \int_{0}^{\lambda}\left( \frac{e^{(w_1+w_2-a)t}}{(e^{w_1t}-1)(e^{w_2t}-1)}-\frac{1}{w_1w_2t^2}-\frac{-2a+w_1+w_2}{2w_1w_2t}\right)t^{s-1}dt  \\
		&\quad +\frac{1}{w_1w_2}\frac{\lambda^{s-2}}{s-2} + \frac{-2a+w_1+w_2}{2w_1w_2}\frac{\lambda^{s-1}}{s-1}.
	\end{split}
	\end{eqnarray*}
\end{lemm}
\begin{pf}
\begin{rm}
	The proof is similar to that of {\rm Lemma \ref{inte1}}. From {\rm (\ref{積分表示})}, we have
	\begin{eqnarray*}
	\begin{split}
		\Gamma(s)\zeta_2(s,\bm{w},a)&= \int_{\lambda}^{\infty} \frac{e^{(w_1+w_2-a)t}}{(e^{w_1t}-1)(e^{w_2t}-1)}t^{s-1}dt \\
		&\quad + \int_{0}^{\lambda}\left( \frac{e^{(w_1+w_2-a)t}}{(e^{w_1t}-1)(e^{w_2t}-1)}-\frac{1}{w_1w_2t^2}-\frac{-2a+w_1+w_2}{2w_1w_2t}\right)t^{s-1}dt  \\
		&\quad +\int_{0}^{\lambda} \left(\frac{1}{w_1w_2t^2}+\frac{-2a+w_1+w_2}{2w_1w_2t}\right)t^{s-1}dt. \\
	\end{split}
	\end{eqnarray*}
Then the first integral is holomorphic in the whole complex plane, the second integral is holomorphic in the half-plane $\sigma > 0$, and the third integral is rational, that is, $\Gamma(s)\zeta_2(s,\bm{w},a)$ is meromorphically continued to $\sigma > 0$.
Thus we obtain this lemma.
\qed\end{rm}
\end{pf}
\begin{lemm}\label{fwnlem} Let  $n$ be a positive integer. We define $f_{\bm{w},n}(a)$ by
	\begin{equation*}\label{f}
		f_{\bm{w},n}(a):=(w_1+w_2)^2a^n-w_1^2(a-w_2)^n-w_2^2(a-w_1)^n,~~~~~a > 0.
	\end{equation*}
	Then we have 
	\begin{enumerate}[{\rm (i)}]
		\item $f_{\bm{w},n}\left(\frac{w_1+w_2}{2}\right) >0.$
		\item $f_{\bm{w},n}(a)$ is monotone increasing for $a \ge (w_1+w_2)/2$.
	\end{enumerate}
Therefore, from {\rm (i)} and {\rm (ii)}, $f_{\bm{w},n}(a)>0$ for $a \ge (w_1+w_2)/2$.
\end{lemm}
\begin{pf} 
\begin{rm}
{\rm (i)} Obviously, one has
\begin{equation*}\label{fwn}
	f_{\bm{w},n}\left(\frac{w_1+w_2}{2}\right)=\frac{1}{2^n}\left((w_1+w_2)^{n+2}-w_1^2(w_1-w_2)^n-w_2^2(w_2-w_1)^n \right).
\end{equation*}
There exists a constant $c >0$ such that
\begin{equation*}
	w_2=cw_1.
\end{equation*}
Then we have 
\begin{equation*}
	f_{\bm{w},n}\left(\frac{w_1+cw_1}{2}\right)=\frac{w_1^{n+2}}{2^n}\left((1+c)^{n+2}-(1-c)^n-c^2(c-1)^n\right).
\end{equation*}
For $n \ge 1$, we have 
\begin{eqnarray*}
	&&w_1^{n+2}/2^n >0 \\
\end{eqnarray*}
and
\begin{eqnarray*}
	&&(1+c)^{n+2}-(1-c)^n-c^2(c-1)^n \\
	&&\geq (1+c)^{n+2} - (1+c^2)|c-1|^n \\
	&& > (1+c)^{n+2}-(1+2c+c^2)(1+c)^n \\
	&& = (1+c)^{n+2} -  (1+c)^{n+2}= 0.
\end{eqnarray*}
From these inequalities, $f_{\bm{w},n}\left((w_1+w_2)/2\right)>0$.

(ii) We prove (ii) of this lemma by induction on $n \ge 1$. Suppose $n=1$. We immediately see that 
$f_{\bm{w},1}(a)=w_1w_2(2a+w_1+w_2)$, which is monotone increasing for $a \ge (w_1+w_2)/2$.
Hence we assume that the case of $n=k$ holds and consider the case of $n=k+1$. 
The first derivative of $f_{\bm{w}, k+1}$ is given as
\begin{equation}
	\frac{d}{da}f_{\bm{w},k+1}(a) = (k+1)f_{\bm{w},k}(a).
\end{equation}
Then, using (i) of Lemma \ref{fwnlem} and the assumption of the induction, we obtain 
\begin{equation}
	\frac{d}{da}f_{\bm{w},k+1}(a) >0,~~~~~a \ge (w_1+w_2)/2.
\end{equation}
Therefore we find that $f_{\bm{w},k+1}(a)$ is also monotone increasing for $a \ge (w_1+w_2)/2$.
Thus we obtain the conclusion.

From (i) and (ii) we obtain $f_{\bm{w},n}(a)>0$ for $a \ge (w_1+w_2)/2$.
\qed\end{rm}
\end{pf}

\begin{lemm}\label{g}
If $a \ge  (w_1+w_2)/2$, then the function $G_2(a, \bm{w}, t)$ is negative for all $t>0$.
\end{lemm}

\begin{pf}
\begin{rm}
	Suppose $a \ge (w_1+w_2)/2$. Obviously, we have $w_1w_2t^2(e^{w_1t}-1)(e^{w_2t}-1)>0$ for all $t > 0$. Thus we have only to consider 
	\begin{equation*}
		g_2(a, \bm{w}, t) := w_1w_2t^2(e^{w_1t}-1)(e^{w_2t}-1)G_2(a, \bm{w}, t).
	\end{equation*}
	By differentiating with respect to $t$, we have
	\begin{eqnarray*}
		g_2'(a, \bm{w}, t) &=& 2w_1w_2te^{(w_1+w_2-a)t}+w_1w_2(w_1+w_2-a)t^2e^{(w_1+w_2-a)t} \\
		&&-(w_1+w_2)e^{(w_1+w_2)t}+w_1e^{w_1t}+w_2e^{w_2t}, \\
		g_2''(a,\bm{w},t) &=& 2w_1w_2e^{(w_1+w_2-a)t}+4w_1w_2(w_1+w_2-a)te^{(w_1+w_2-a)t} \nonumber\\
		&&+w_1w_2(w_1+w_2-a)^2t^2e^{(w_1+w_2-a)t} \nonumber \\
		&&-(w_1+w_2)^2e^{(w_1+w_2)t}+w_1^2e^{w_1t}+w_2^2e^{w_2t}.
	\end{eqnarray*} 
	Obviously, one has 
	\begin{equation}
	 g_2(a, \bm{w}, 0) = g'_2(a, \bm{w}, 0) = g''_2(a, \bm{w}, 0) = 0.
	 \end{equation}
Now we consider the function
	\begin{eqnarray}
		&&e^{-(w_1+w_2-a)t}g_2''(a,\bm{w},t) \nonumber \\
		&=& 2w_1w_2+4w_1w_2(w_1+w_2-a)t +w_1w_2(w_1+w_2-a)^2t^2\nonumber\\
		&&-(w_1+w_2)^2e^{at}+w_1^2e^{(a-w_2)t}+w_2^2e^{(a-w_1)t}.
	\end{eqnarray}
	By the Taylor expansion of $e^x$, one has
	\begin{eqnarray*}
		&&(w_1+w_2)^2e^{at}-w_1^2e^{(a-w_2)t}-w_2^2e^{(a-w_1)t} \nonumber \\
		&=&\sum_{n=0}^{\infty}\left((w_1+w_2)^2a^n-w_1^2(a-w_2)^n-w_2^2(a-w_1)^n\right)\frac{t^n}{n!} \nonumber \\
		&=&2w_1w_2+w_1w_2(2a+w_1+w_2)t+w_1w_2 \left(a^2+(w_1+w_2)a-w_1w_2\right)t^2 \nonumber \\
		&&+\sum_{n=3}^{\infty}\left((w_1+w_2)^2a^n-w_1^2(a-w_2)^n-w_2^2(a-w_1)^n\right)\frac{t^n}{n!}.
	\end{eqnarray*}
	Hence we have
	\begin{eqnarray*}
	&&e^{-(w_1+w_2-a)t}g_2''(a,\bm{w},t) \nonumber \\
		&=& -w_1w_2(6a-3(w_1+w_2))t -w_1w_2(3a(w_1+w_2)-w_1w_2-(w_1+w_2)^2)t^2\nonumber\\
		&&-\sum_{n=3}^{\infty}\left((w_1+w_2)^2a^n-w_1^2(a-w_2)^n-w_2^2(a-w_1)^n\right)\frac{t^n}{n!}.
	\end{eqnarray*}
Let $a \ge (w_1+w_2)/2$. Then 
	\begin{eqnarray*}
		&&w_1w_2(6a-3(w_1+w_2))  \geq 0 ,\\
		&& w_1w_2(3a(w_1+w_2)-w_1w_2-(w_1+w_2)^2) \geq 0.
	\end{eqnarray*}
By Lemma \ref{fwnlem}, 
\begin{equation*}\label{eq}
	(w_1+w_2)^2a^n-w_1^2(a-w_2)^n-w_2^2(a-w_1)^n > 0,~~~~~n \ge 3.
\end{equation*}
From these inequalities, $e^{-(w_1+w_2-a)t}g_2''(a,\bm{w},t)<0$ for all $t>0$. Thus we obtain $g_2(a,\bm{w},t)<0$ and $G_2(a, \bm{w}, t) <0$ for all $t > 0$.
\qed\end{rm}
\end{pf}

\begin{proof}[{\bf Proof of (1) of Theorem \ref{main}}]
\begin{rm}
	Let $0 < a < (w_1+w_2)/2$. By {\rm Lemma \ref{inte2}}, we see that
		\begin{eqnarray}\label{limit}
			\lim_{\sigma \rightarrow 2-}\zeta_2(\sigma, \bm{w}, a) = -\infty , \lim_{\sigma \rightarrow 1+}\zeta_2(\sigma, \bm{w}, a) = \infty.
		\end{eqnarray}
	Thus $\zeta_2(\sigma, \bm{w}, a)$ has zeros in $(1,2)$.
	
	Secondly suppose $a \ge (w_1+ w_2)/2$. Then we have 
\[
	\Gamma(\sigma)\zeta_2(\sigma, \bm{w}, a) = \int_{0}^{\infty} G_2(a, \bm{w}, t)t^{\sigma-1}dt,~~~1<\sigma<2
\]
	by the integral expression (\ref{inte}). It is well-known that $\Gamma(\sigma) > 0$ for any $1 < \sigma < 2$. Thus we obtain 
		\begin{equation*}
			\zeta_2(\sigma, \bm{w}, a)<0
		\end{equation*}
	for all $1<\sigma<2$ by Lemma \ref{g} and the integral expression above. Therefore $\zeta_2(\sigma, \bm{w}, a)$ does not vanish in the interval $(1, 2)$ when $a \ge (w_1+w_2)/2$.
\end{rm}
\end{proof} 

Next, we show (2) of Theorem \ref{main} by using the method of Endo-Suzuki (see \cite[(1) of Theorem 1.2]{es}).

\begin{lemm}\label{Gsimple}
	For any $0<a<(w_1+w_2)/2$, there exists a positive $t_0$ such that 
	\begin{eqnarray*}
		G_2(a, \bm{w}, t_0) = 0,~~G_2(a, \bm{w}, t) > 0 ~(for~0<t<t_0), ~~G_2(a, \bm{w}, t) < 0 ~(for~t_0 < t).
	\end{eqnarray*}
\end{lemm}
\begin{pf}
\begin{rm}
	We consider $g_2(a, \bm{w}, t)$ in the proof of Lemma \ref{g}. We recall that
	\begin{eqnarray*}
		g_2'(a, \bm{w}, t) &=& 2w_1w_2te^{(w_1+w_2-a)t}+w_1w_2(w_1+w_2-a)t^2e^{(w_1+w_2-a)t} \\
		&&-(w_1+w_2)e^{(w_1+w_2)t}+w_1e^{w_1t}+w_2e^{w_2t}.
	\end{eqnarray*}
Put
	\begin{eqnarray*}
		h_2(a, \bm{w}, t) &:=& e^{(a-w_1-w_2)t}g_2'(a, \bm{w}, t) \\
		&=& 2w_1w_2t+w_1w_2(w_1+w_2-a)t^2 \\
		&&-(w_1+w_2)e^{at}+w_1e^{(a-w_2)t}+w_2e^{(a-w_1)t}.
	\end{eqnarray*}
By differentiating $h_2(a, \bm{w}, t)$ with respect to $t$, we have
	\begin{eqnarray*}
		h_2'(a, \bm{w}, t) &=& 2w_1w_2+2w_1w_2(w_1+w_2-a)t \\
		&&-a(w_1+w_2)e^{at}+w_1(a-w_2)e^{(a-w_2)t}+w_2(a-w_1)e^{(a-w_1)t}, \\
		h_2''(a, \bm{w}, t) &=& 2w_1w_2(w_1+w_2-a) \\
		&&-a^2(w_1+w_2)e^{at}+w_1(a-w_2)^2e^{(a-w_2)t}+w_2(a-w_1)^2e^{(a-w_1)t},\\
		h_2'''(a, \bm{w}, t) &=& -a^3(w_1+w_2)e^{at}+w_1(a-w_2)^3e^{(a-w_2)t}+w_2(a-w_1)^3e^{(a-w_1)t} .
	\end{eqnarray*}
Then we obtain that
	\begin{eqnarray*}
		h_2(a, \bm{w}, 0) = 0,&&\lim_{t \rightarrow \infty} h_2(a, \bm{w}, t) = -\infty, \\
		h_2'(a, \bm{w}, 0) = 0,&&\lim_{t \rightarrow \infty} h_2'(a, \bm{w}, t) = -\infty ,
	\end{eqnarray*}
	\begin{eqnarray*}
		h_2''(a, \bm{w}, 0) = -6w_1w_2\left(a-\frac{w_1+w_2}{2} \right) > 0,&&\lim_{t \rightarrow \infty} h_2''(a, \bm{w}, t) = -\infty, \\
	\end{eqnarray*}
	\begin{eqnarray*}
		h_2'''(a, \bm{w}, t) &=& -a^3(w_1+w_2)e^{at}+w_1(a-w_2)^3e^{(a-w_2)t}+w_2(a-w_1)^3e^{(a-w_1)t} \\
		&<& -a^3(w_1+w_2)e^{at}+w_1a^3e^{(a-w_2)t}+w_2a^3e^{(a-w_1)t} \\
		&<& -a^3(w_1+w_2)e^{at}+w_1a^3e^{at}+w_2a^3e^{at} = 0.
	\end{eqnarray*}
Hence we have the conclusion.
\end{rm}
\qed\end{pf}

\begin{lemm}\label{strictly}
	Let $0 < a < (w_1+w_2)/2$ and $t_0 > 0$ as in Lemma \ref{Gsimple}. Then the function
	\begin{equation*}
		t_0^{-\sigma} \Gamma(\sigma)\zeta_2(\sigma, \bm{w}, a) = \int_{0}^{\infty} G_2(a, \bm{w}, t) \left(\frac{t}{t_0}\right)^\sigma \frac{dt}{t}
	\end{equation*}
	is strictly decreasing for $1 < \sigma < 2$.
\end{lemm}
\begin{pf}
\begin{rm}
	See \cite[Lemma 2.3]{es}.
\end{rm}
\end{pf}

\begin{proof}[{\bf Proof of (2) of Theorem \ref{main}}]
Let $0<a<(w_1+w_2)/2$. Put 
\begin{equation*}
	F_2(\sigma, \bm{w}, a) = t_0^{-\sigma}\Gamma(\sigma)\zeta_2(\sigma, \bm{w}, a).
\end{equation*}
Using (\ref{limit}), we have
\begin{eqnarray*}
	\lim_{\sigma \rightarrow 1+} F_2(\sigma, \bm{w}, a) = \infty, \lim_{\sigma \rightarrow 2-} F_2(\sigma, \bm{w}, a) = -\infty.
\end{eqnarray*}
By Lemma \ref{strictly}, $F_2(\sigma, \bm{w}, a)$ has one zero in the interval $(1, 2)$. Let $\beta^{(2)}(\bm{w}, a)$ denote this unique zero contained in $(1, 2)$. Since $t_0^{-\sigma} \Gamma(\sigma)$ has no zero and no pole in $(1, 2)$, the zeros of $F_2(\sigma, \bm{w}, a)$ and $\zeta_2(\sigma, \bm{w}, a)$ and their orders coincide. Hence it suffices to show that $\beta^{(2)}(\bm{w}, a)$ is a simple zero of $F_2(\sigma, \bm{w}, a)$. For any $0 < a < (w_1+w_2)/2$, we have
	\begin{equation*}
		\frac{\partial}{\partial \sigma} F_2(\sigma, \bm{w}, a) = \left(\int_{0}^{t_0} + \int_{t_0}^{\infty} \right)G_2(a, \bm{w}, t) \left(\frac{t}{t_0}\right)^\sigma \log \left(\frac{t}{t_0}\right) \frac{dt}{t} < 0
	\end{equation*}
since, for any $t \neq t_0$, Lemma \ref{Gsimple} implies
	\begin{equation*}
		G_2(a, \bm{w}, t)\left(\frac{t}{t_0}\right)^\sigma \log \left(\frac{t}{t_0}\right) < 0.
	\end{equation*}
This completes the proof.
\end{proof}

\section{Some problems}
We conclude this paper with the discussion of some further problems.
\subsection{Real zeros of the Barnes $r$-ple zeta function}
First, we define $_N G_r(x, \bm{w}, t)$ by
\begin{defi}Let $N \ge -r$ be an integer. We set
\begin{equation}
	_NG_r(x,\bm{w},t) := \frac{e^{(w_1+\cdot\cdot\cdot + w_r-x)t}}{(e^{w_1t}-1)\cdot\cdot\cdot (e^{w_rt}-1)}-\sum_{k=0}^{N+r}B_k^{(r)}(x,\bm{w}) \frac{t^{k-r}}{k!}.
\end{equation}
\end{defi}
Then we can show the following lemmas.
\begin{lemm} Let $N \ge -r$ be an integer. Then we have
	\begin{equation}
		\Gamma(s)\zeta_r(s, \bm{w}, a) = \int_{0}^{\infty} {} _NG_r(a,\bm{w},t) t^{s-1}dt,~~~-N-1<\Re(s)<-N.
	\end{equation}
\end{lemm}

\begin{lemm}\label{一般} Let $N \ge -r$ be an integer and fix a $\lambda \in \mathbb{R}$ which satisfies $0< \lambda < \min \left\{2\pi/w_1, ... , 2\pi/w_r \right\}$. Then we have
	\begin{eqnarray*}
		\Gamma(s)\zeta_r(s, \bm{w}, a) &=& \int_{\lambda}^{\infty} \frac{e^{(w_1+\cdot\cdot\cdot +w_r-a)t}}{(e^{w_1t}-1)\cdot\cdot\cdot (e^{w_rt}-1)}t^{s-1}dt \\
		&&+  \int_{0}^{\lambda} {} _{N+1}G_r(a,\bm{w},t) t^{s-1}dt  \\
		&& + \sum_{k=0}^{N+r+1} \frac{1}{k!} B_k^{(r)}(a,\bm{w}) \frac{\lambda^{s+k-r}}{s+k-r},~~~\Re(s)>-N-2.
	\end{eqnarray*}
\end{lemm}


Finally, we obtain a weak extension of (1) of Theorem \ref{main} from Lemma \ref{一般}.
\begin{theo}\label{拡張} Let $N \ge -r$ be an integer. Then we have the following:

		If $B_{r+N+1}^{(r)}(a, \bm{w})B_{r+N}^{(r)}(a, \bm{w})>0$, then $\zeta_r(\sigma, \bm{w}, a)$ has zeros in $(-N-1, -N)$.
\end{theo}
\begin{pf}
\begin{rm}
	We can show Theorem \ref{拡張} from Lemma \ref{一般} and the intermediate value theorem.\\
\end{rm}
\qed\end{pf}
When $r = 1$, this theorem agrees with Theorem \ref{m1}.

Thus, we can conjecture the following.
\begin{con} Let $N \ge -r$, be an integer. Then $\zeta_r(\sigma, \bm{w}, a)$ has zeros in $(-N-1,-N)$ if and only if 
$B_{r+N+1}^{(r)}(a, \bm{w})B_{r+N}^{(r)}(a, \bm{w})>0$.
\end{con}

\begin{ack}
\
{\rm The author is deeply grateful to Prof. Kohji Matsumoto, Prof. Hirofumi Tsumura, Prof. Takashi Nakamura, Mr. Toshiki Matsusaka, Mr. Yuta Suzuki and Mr. Kenta Endo for their helpful comments and interest.}
\end{ack}

GRADUATE SCHOOL OF MATHEMATICS, NAGOYA UNIVERSITY, FUROCHO, CHIKUSAKU, NAGOYA 464-8602, JAPAN

\textit{E-mail address}: d20002i@math.nagoya-u.ac.jp
\end{document}